\documentclass{ifacconf}
\usepackage{graphicx}     
\usepackage{natbib}  
\usepackage{color,theorem}

\usepackage{amssymb}
\usepackage{amsmath}
\usepackage{mathrsfs}
\usepackage{algorithm}
\usepackage{algorithmic}
\usepackage{bm}
\usepackage{xcolor}
\usepackage{color}
\usepackage{pgfplots}
\usepackage{theorem}
\usepackage{enumitem}

\makeatletter
\let\old@ssect\@ssect % Store how ifacconf defines \@ssect
\makeatother

 \usepackage{hyperref}
 
\makeatletter
\def\@ssect#1#2#3#4#5#6{%
  \NR@gettitle{#6}% Insert key \nameref title grab
  \old@ssect{#1}{#2}{#3}{#4}{#5}{#6}% Restore ifacconf's \@ssect
}
\makeatother 
 
\setitemize[1]{itemsep=5pt,partopsep=3pt,parsep=\parskip,topsep=5pt}

\begin{document}
\begin{frontmatter}

\title{Teaching MPC: \\ Which Way to the Promised Land? } 

\author[First]{Timm Faulwasser} 
\author[Second]{Sergio Lucia}
\author[Third]{Moritz Schulze Darup} 
\author[Fourth]{Martin Mönnigmann}

\address[First]{Institute of Energy Systems, Energy Efficiency and Energy Economics,  TU Dortmund University, Dortmund, Germany \\	{\tt {timm.faulwasser}@ieee.org}}
\address[Second]{Laboratory of Process Automation Systems,  TU Dortmund University, Dortmund, Germany \\	{\tt {sergio.lucia}@tu-dortmund.de}}
\address[Third]{Control and Cyberphysical Systems, TU Dortmund University, Dortmund, Germany 	\\ {\tt {moritz.schulzedarup}@tu-dortmund.de}}
\address[Fourth]{Automatic Control and Systems Theory, Ruhr University Bochum, Bochum, Germany \\ {\tt {martin.moennigmann}@rub.de}}
%\thanks[footnoteinfo]{Th}

\begin{abstract}% Abstract of not more than 250 words.
Since the earliest conceptualizations by Lee and Markus, and Propoi in the 1960s, Model Predictive Control (MPC) has become a major success story of systems and control with respect to industrial impact \textit{and} with respect to continued and wide-spread research interest. The field has evolved from conceptually simple linear-quadratic (convex) settings in discrete and continuous time to nonlinear and distributed settings including hybrid, stochastic, and infinite-dimensional systems. Put differently, essentially the entire spectrum of dynamic systems can be considered in the MPC framework with respect to both---system theoretic analysis and tailored numerics. Moreover, recent developments in machine learning also leverage MPC concepts and learning-based and data-driven MPC have become highly active research areas.  

However, this evident and continued success renders it increasingly complex to live up to industrial expectations while enabling graduate students for state-of-the-art research in teaching MPC.
Hence, this position paper attempts to trigger a discussion on teaching MPC.
To lay the basis for a fruitful debate, we subsequently investigate the prospect of covering MPC in undergraduate courses; we comment on teaching textbooks; and we discuss the increasing complexity of research-oriented graduate teaching of~MPC.
\end{abstract}

\begin{keyword}
Model Predictive Control, Teaching
\end{keyword}

\end{frontmatter}

\section{Introduction}
%\typeout{col width is \the\columnwidth}

The idea of Model Predictive Control (MPC)---or receding-horizon control---can at least 
be traced back to  \cite{Propoi63} and to \cite{Lee67}. Roughly 50-60 years later, its success in industrial applications cannot be ignored. 
Among other features, MPC stands out due to the effective combination of the three following key items: 

\begin{itemize}[leftmargin=6.8mm] \setlength{\itemsep}{0pt}
\item[(a)] simplicity of the conceptual idea, 
\item[(b)] efficacy of handling complex systems, and
\item[(c)] flexibility to dock with emerging trends in systems and control, general engineering, and computer science. 
\end{itemize}

With respect to item (a), many control  educators will agree that the conceptual idea of MPC is rather straightforwardly conveyed to students---both %either
undergraduate and graduate ones. Similar experiences can be made in communicating to industrial practitioners not previously accustomed to MPC. Moreover, it stands to reason that the mathematical prerequisites required to understand a linear-quadratic MPC controller are not more complex than those needed for 
for frequency-domain techniques in classical PID control. Yet, simplicity alone does not guarantee success. Hence, the second important aspect of MPC is that it has proven itself to be a very useful method in manifold applications, cf. item (b). 
While the classical optimal control route to deriving optimal \textit{feedbacks} via  the Hamilton-Jacobi-Bellman Equation or via the Pontryagin Maximum Principle  is only viable in quite specific settings---i.e. mostly linear dynamics or rather low-dimensional nonlinear systems---the concept of MPC yields receding-horizon feedback for a large class of systems.  

Long gone are times when  real-time feasible MPC and NMPC applications were limited to slow process systems, see \cite{Qin00,Qin03a} for overviews of early industrial applications. Without any claim of completeness, milestones in terms of computational speed can be marked by  \cite{Diehl01a}, who demonstrated NMPC on a lab-scale system with sampling time 10s, \cite{Jerez14a} who report MHz sampling rates for MPC and \cite{Monnigmann2010b} (who target GHz sampling rates using explicit MPC). Moreover, \cite{Houska11a} report kHz sampling rates for nonlinear systems. Enablers for this progress have been papers such as   \citep{Bemporad02a,Seron2003,Wang2010}, while \cite{Ohtsuka02a} proposed code generation for NMPC. 
Hence, it is far from surprising that nowadays, powerful open-source and commercial software tools for MPC (see e.g. \cite{Lofberg04,Houska11a,MPT3,MatlabMpcToolbox,Andersson19a})\footnote{The authors do not claim this list to be comprehensive. See, e.g., Table 1 in \citep{Findeisen2018} for an overview of tools and the corresponding timeline.} enable readily solving problems of a complexity that would have been publication worthy not too long ago. 

\begin{figure}
\begin{center}
\includegraphics[width=6.881cm]{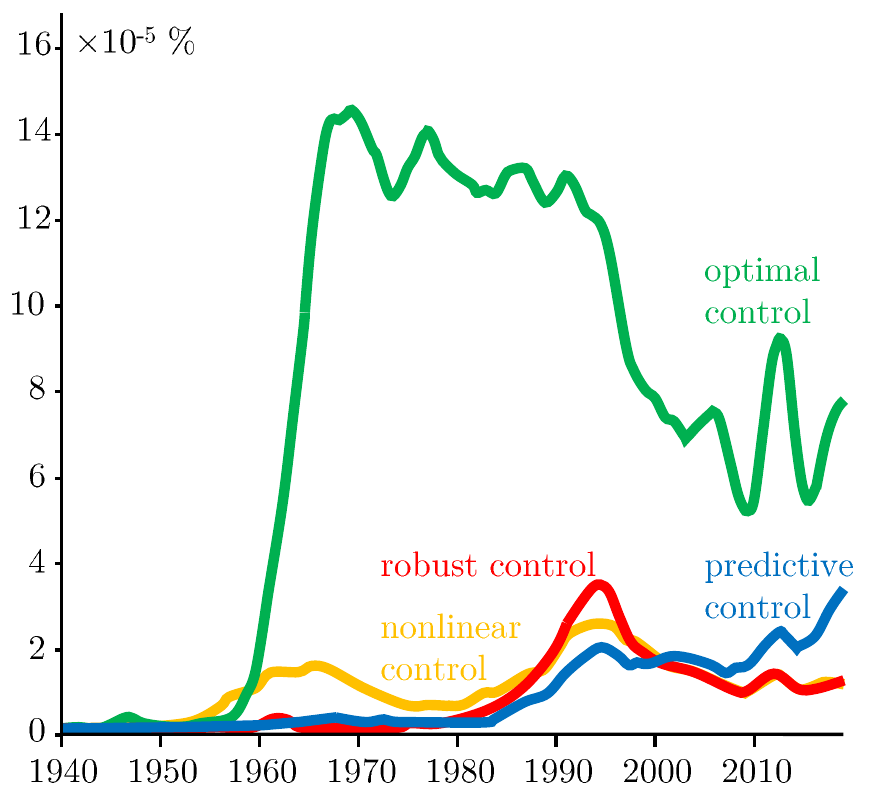}
\end{center}
\caption{Frequency of search terms in English books according to Google Ngram Viewer. Data is normalized by the number of books published in each year.}\label{fig:ngram}
\end{figure}

Turning to item (c), it is interesting to observe that while most research topics are subject to activity cycles---the infamous AI-winter \citep{AIwinter} being a prominent example---MPC has not seen such cycles yet. 
As evident from Fig.~\ref{fig:ngram}, closely related topics such as \textit{optimal control} and \textit{dynamic programming} have seen periods of varying scientific interest (as measured by the Google Ngram viewer, which counts frequency of search terms in English texts). MPC, in contrast, has seen a steady growth since the early 1980s.
It is fair to ask why  a research field in systems and control seemingly grows void of apparent activity cycles. The answer can be found in the flexibility of  MPC to dock with other control trends, which in turn led to the emergence of various branches such as $\{$nonlinear, stochastic and robust, hybrid, distributed, economic, data-driven, learning$\}$ MPC. Likewise,  MPC finds use in a variety of application domains: process engineering, mechatronics and robotics, aeronautics, logistics, smart grids and energy systems, finance, cyber-security, etc.  Arguably, this twofold diversity is a unique advantage of MPC. 

While this flexibility has catalyzed 
research on MPC for 
more than three decades---one might go as far as saying that  MPC  is adaptive with respect to parallel research trends---it also renders teaching MPC increasingly complex, though necessary. This complexity does not mean that MPC courses should be postponed to the graduate level. The necessity to teach MPC
stems on the one hand from the requirement that higher education in engineering (undergraduate and graduate alike) should prepare for careers in industry. Systems and control would fail in this regard if MPC was not part of  the core curriculum in control. On the other hand, it becomes increasingly difficult on the graduate level to adhere to \textit{Humboldt's ideal of combining  research and education}
%\footnote{We refer to  \cite{Pritchard04a} for more details on the long term impact of the Humboldtian ideal on modern German universities.}  
and to cover the entire spectrum of MPC branches. 
At the same time, hardly any research group will address the complete spectrum of MPC topics in their research. This raises the questions of what 
the research topics are that can be covered in a graduate course, and what their relation with topics in systems and control and other neighboring disciplines are.

It is \textit{not} the purpose of this paper
to provide the receding-horizon optimal solution to teaching MPC. Neither does it elaborate the question of whether MPC should be taught in dedicated courses. Rather, it asks where MPC positions itself with respect to the undergraduate-graduate teaching divide. Specifically, the paper intends to trigger discussion within the MPC community on two aspects: 
\begin{itemize}\setlength{\itemsep}{0pt}
\item[(i)] prospect and limits of  teaching MPC to {\it undergraduate} engineering students, and
\item[(ii)] synergies of advanced MPC courses with neighboring control topics and with disciplines beyond the engineering realm (such as computer science or applied mathematics).
\end{itemize}
Both aspects call for a brief discussion of existing textbooks on MPC. Hence, Section~\ref{subsec:Textbooks} gives a short summary of the books the authors have relied on in preparing their undergraduate and graduate courses so far.   

The remainder of this paper is structured as follows:
Section \ref{sec:LQMPC} proposes a curriculum for an undergraduate curriculum, including a brief discussion of existing textbooks in Section~\ref{subsec:Textbooks}. 
Section~\ref{sec:NMPC} lists advanced MPC topics that may be covered in graduate courses. 
The paper closes with hypotheses on teaching MPC in Section \ref{sec:hypo}, which are meant to serve as the basis of a plenary discussion at the 
7$^{\text{th}}$ IFAC Conference on Nonlinear Model Predictive Control 2021.

\section{The Prospects of Undergraduate MPC Teaching} \label{sec:LQMPC}

As pointed out, MPC is a powerful, widely applied, and conceptually easy control method. Thus, it stands to reason that it should play a prominent role in control education. Hence, the question arises whether it can be addressed early in the curriculum in order to leave enough time for advanced control courses. 
However, while conceptually easy, linear-quadratic MPC combines elements from distinct  fields such as linear algebra, systems theory, and numerical optimization. For effectively teaching basic MPC to undergraduates, it is thus important to determine the bare necessities regarding students' prior knowledge.

Clearly, 
the assessment requires additional information about the underlying undergraduate curriculum. 
The situation will likely differ among institutions and even among different departments in the same institution.  
This holds in particular since MPC has attracted sustained interest from engineering, mathematics and, more lately, computer science departments.
We believe, however, that  certain pivotal insights are independent of specific situations and we attempt to summarize them below. Such a summary is always based on the combined experience of the authors and, thus, never truly 
objective.  Whether these insights form indeed a viable nucleus of undergraduate MPC teaching and how certain bottlenecks should be tackled deserves further discussion within the MPC community, respectively, it is for the reader to decide. 

\subsection{Bare necessities of prior knowledge}\label{subsec:preknowledge}

There are only a handful of requirements for an undergraduate MPC course. Students should be familiar with 
\begin{itemize}\setlength{\itemsep}{0pt}
\item basics of linear algebra (matrix computations, eigenvalues, vector spaces), 
\item analysis in $\mathbb{R}^n$ (gradients, Hessians, extrema), %of functions), 
\item basic linear systems theory (state space models, stability of LTI systems), and
\item control fundamentals (feedback vs. feed-forward). %controls). 
\end{itemize} 
Put differently, an undergraduate MPC course should not be responsible for laying the mathematical foundations or for providing first insights into systems and control. 
One may even say that students must have learned about the pivotal nature of feedback in order to appreciate optimizing over a long horizon and dispensing with the largest part of the open-loop predictions.
Clearly, there is a long and continuously growing  wish list of additional prior knowledge beneficial for an MPC course such as first contacts with basic LTI optimal control (LQR) or convex optimization (optimality conditions). However, we believe that these essentials can also reasonably covered by the MPC course itself.

\subsection{Modeling and system identification}

Although the importance of (physical or data-driven) state-space models for MPC is undeniable, their derivation is often skipped when teaching MPC. This may turn out to be critical in undergraduate courses because of two problems: (i) Students with no further background in modeling of dynamic systems can underestimate the challenge of modeling for MPC and other control methods. 
(ii) They could also be unlikely to appreciate a model-based method if the origin of the models remains unclear. One may argue, however, it is enough to show by example that systematic methods for deriving models adequate for MPC exist, and to encourage students to attend an additional course that covers modeling in greater depths. To this end, it suffices to revisit modeling based on conservation laws (conservation of mass, momentum and energy, Kirchhoff's laws depending on the preferences of the lecturer), the linearization of dynamical systems, and linear and nonlinear least-squares optimization, which can be used for parameter estimation or identification of the system matrices. These topics need not be covered in great theoretical depths, but a pragmatic introduction, ideally combined with a hands-on exercise with numerical tools, is sufficient to demonstrate that models of the required type can be derived from data or from first principles combined with data.

\subsection{Unconstrained and constrained optimization}

The natural and straightforward inclusion of state and input constraints as well as the built-in applicability to MIMO systems undoubtedly are key features of MPC. 
In this context, it is important to point out that solving a control task while explicitly taking constraints into account is fundamentally different from the common approach of augmenting %patching 
``unconstrained'' controllers with subsequent modifications (such as, e.g., PID with anti-windup).\footnote{This is not to say that modern anti wind-up schemes do not have their merit and place in industrial applications, cf. \citep{Galeani09a}. Yet, such schemes are also non-trivial to analyze thoroughly and thus involve educational challenges of their own. A detailed discussion is, however, beyond the scope of the present paper.}
Analogously, under the assumption that students have no background on convex optimization, it is essential to work out the fundamental differences between unconstrained and constrained optimization. 
While the latter is crucially important for modern MPC,\footnote{Indeed, traditional unconstrained MPC variants such as GPC may or may not be included in an undergraduate course.} the former often provides useful links to prior knowledge such as the analysis of extreme values or LQR.

Apart from these general insights, the level of detail regarding constrained convex optimization and constrained quadratic programming will depend on the focus of the course itself and subsequent courses. The spectrum might range from the application of  toolboxes---see, e.g.,  \cite{Lofberg04,Houska11b,MPT3,MatlabMpcToolbox,Andersson19a}---over
problem modelling and formulation combined with black-box QP solvers, to including the coding of low-complexity QP-solvers---e.g. Nesterov's fast gradient method; see \cite{Nesterov1983} and \cite{Richter09,ifat:zometa12} for its use in MPC)---as practical (flipped classroom or exercises) course elements. Practical coding efforts will benefit from 
 the availability of school-wide licenses of computing environments such as Matlab and previous introduction  to numerical computations. 
However, as costly school-wide Matlab licensing schemes may not be available at all institutions, one can also straightforwardly  resort to open-source alternatives such as Python or Octave.
 In any case, there is significant prospect for an undergraduate amalgam of MPC teaching \textit{and} an introduction to numerical computations.  Since at many institutions the vast majority of undergraduate students will likely have access to private laptop computers, coding tasks can easily be blended into the lectures; thus implementing the concept of active learning~\citep{Hackathorn11}.

\subsection{Stability versus optimality}

Another crucial observation for MPC is that optimality only holds with respect to the formulated mathematical problem.
Receding-horizon optimization does not necessarily guarantee typical design criteria such as stability or robustness, see \cite{Mayne00a}.
This can be easily illustrated with suitable numerical examples many of which are readily available. As it is well-known, designing the MPC problem such that stability is guaranteed is significantly more difficult than merely stating a plausible optimization problem.
Nevertheless, it is conceptually easy to motivate the two most-common approaches---the inclusion of a ``safe'' terminal set or the choice of a sufficiently long prediction horizon---as both approaches can easily be illustrated with numerical experiments. 
Whether or not providing formal proofs in an undergraduate course is meaningful, and whether or not proving stability is crucial for educating practitioners-to-be,
are interesting questions that merit discussion.

Another element worth discourse is that not in all applications of MPC closed-loop stability is the intended purpose. In synchronization problems and in scheduling applications---e.g. in energy systems, resource allocation, or finance applications---traditional stability is not aimed for. Yet, while recent developments in economic MPC appear beyond reach for undergraduate teaching, the interplay of stability and optimality has more facets than the well-known Kalman quote on optimality not implying stability~\citep{Kalman60a}.

\subsection{Case studies versus rigorous theory}

Effective undergraduate teaching in engineering departments has to put extra effort into balancing theory and examples. 
 Fortunately, in all engineering domains there exist numerous examples suitable for  undergraduate students. In mechanical engineering, one may consider mobile robots modelled as double integrators or linearized models of quadcopters. In electrical engineering, linearized generator swing equations allow tackling frequency stabilization via linear MPC. Moreover, battery scheduling provides another and easy route to the consideration of problems beyond stabilization. In chemical engineering, control of continuous reactors is a straightforward option. Remarkably, all these examples allow for an easy motivation of constrained MIMO control. 

While the right balance 
of theory and case studies depends on many extrinsic factors---ranging from department culture to teacher preferences and beyond---undergraduate MPC teaching provides the remarkable opportunity to introduce a hands-on, MIMO-applicable, constrained control method without going through all the details of the analysis of MIMO systems. One may go as far as arguing that MIMO LTI control taught from an MPC point of view is easier and more effective than the traditional PID route. Indeed, the latter involves and requires  quite specific concepts such as transmission zeros, control channel decoupling and relative-gain-array analysis. 
In fact, it can be quite insightful for students to uncover the piecewise affine structure of the optimal control law (or, analogously, of the solution to multiparametric QPs).

Moreover, realistic case studies will naturally guide students towards the duality of control problems, i.e., the need to design state estimators. While Kalman filters and observers might shoot with cannons for birds, the classical Luenberger observer does not pose any undergraduate teaching obstacle. Yet, it means that pole placement and state feedback concepts are to be included or to be required.

In summary, we claim that---if balanced nicely between theory and application---teaching MPC on the undergraduate level provides a number of opportunities. Besides the ones mentioned above, and beyond the appeal that this may naturally have for faculty who research MPC, it may also turn out as an attractive student gateway towards graduate studies in  systems and control.

\subsection{Available textbooks}\label{subsec:Textbooks}

Any discussion about teaching MPC is incomplete without reference to available textbooks. Evidently and luckily, the last decade has seen a substantial growth of textbooks on MPC. Moreover, in February 2021, the search string \textrm{predictive control} returned 140--400 results in English on \textrm{amazon.com},\footnote{Depending on the precise settings used.} which are as diverse as re-prints of US patents, books on applications of MPC to specific domains, or research and teaching textbooks. Hence a comprehensive overview of MPC textbooks appears infeasible. Yet, there are undoubtedly established texts and the comparison of which provides indications on the current status. 

Classic and early textbooks on MPC include \cite{Maciejowski2002},~\cite{Rossiter2004} and \cite{Camacho2007}, where the first edition of the latter appeared as early as 1999. Not surprisingly these texts focus mostly on the linear setting. It appears that, in the potentially biased view of the authors,  \cite{Maciejowski2002}  comes quite close to covering all elements needed for an undergraduate MPC teaching endeavor. Yet, at the time of writing it appears to be out of print. Other books that cover most of the elements needed for undergraduate teaching include~\cite{Wang09}, 
\cite{Rawlings2017}, whose first edition appeared in 2009, and~\cite{Borelli2017,Kouvaritakis2016}.
In contrast, \cite{Gruene2017} (first edition in 2013) and \cite{Ellis17a} focus on more advanced topics including economic MPC.

Two facts stand out: there appears to be no exclusively introductory textbook, i.e., no  MPC counterpart of \cite{Astrom21}. This is \textit{not} to be read as a critique of the books mentioned above. Rather all the books referenced---seen in the context of their respective publication dates---are teaching \textit{and} research oriented.

Moreover, when it comes to advanced settings, the books mentioned above, except for \cite{Ellis17a} and \cite{Wang09}, only touch upon continuous-time formulations. The latter aspect is reflective of the fact that MPC research on finite-dimensional systems has developed a bias towards the discrete-time framework.

\begin{figure*}[tp]
\begin{center}
\includegraphics[width=\textwidth]{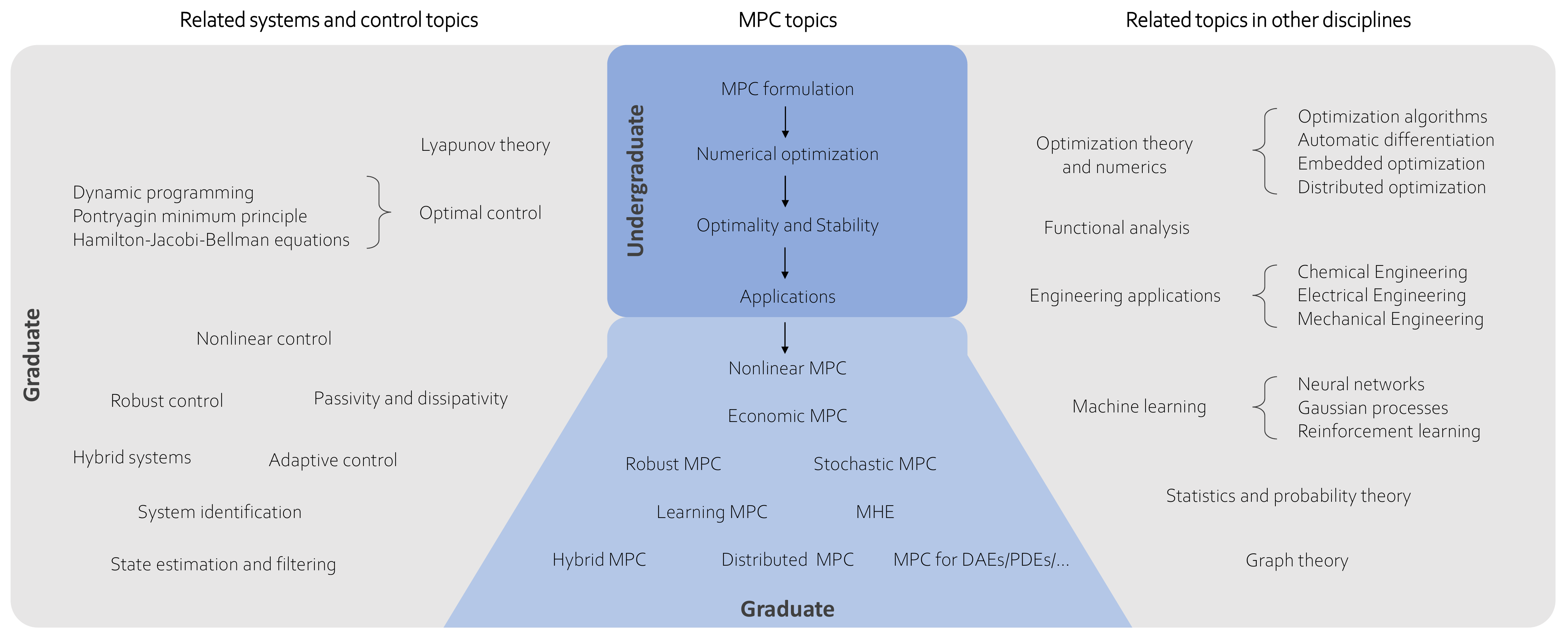}
\end{center}
\caption{Illustration of the main topics covered in a potential MPC lecture for undergraduates and advanced MPC topics for a graduate level course. The possible advanced topics are very broad and often require knowledge from related systems and control fields (left column) as well as neighboring disciplines (right column).}\label{fig:mpc_diagram}
\end{figure*}

\section{Advanced topics in MPC and Synergies with Neighboring Domains} \label{sec:NMPC}

Under the assumptions of the basic prior knowledge described in Section~\ref{subsec:preknowledge}, an MPC course for undergraduate students can be organized in a sequential manner. In this case, the formulation of the MPC problem is followed by notions of optimization, ideas of stability and the demonstration of MPC on simple applications related to different engineering domains, as represented in the central top part of Fig.~\ref{fig:mpc_diagram}.

In clear contrast to this sequential structure, the scope of advanced MPC topics that can be presented in a graduate MPC course is much broader. The wide range of topics leads to many possible pathways and focus areas, which can be established, e.g., depending on the specific engineering department or the research interests of the lecturer. The advanced topics cover usual concepts discussed in most of the MPC textbooks mentioned in Section~\ref{subsec:Textbooks} such as nonlinear MPC including Lyapunov theory, but also more specialized ideas as MPC for hybrid systems or stochastic MPC.

The central bottom part of Fig.~\ref{fig:mpc_diagram} represents the wide variety of advanced MPC topics that can be covered in an MPC graduate course. Observe that this does not only include MPC but also its counterpart, i.e., Moving Horizon Estimation (MHE). It is also interesting to note that the teaching of advanced MPC requires additional concepts from systems and control (left part of Fig.~\ref{fig:mpc_diagram}) as well as from neighboring disciplines such as optimization theory, machine learning or statistics (see right part of Fig.~\ref{fig:mpc_diagram}).

If students have a broader background, a graduate course may safely deviate from  the undergraduate structure. For example, an in-depth discussion of optimal control theory with details on dynamic programming and Pontryagin's minimum principle could be a good starting point to motivate the formulation of the MPC problem. While the central points are probably shared between the graduate and undergraduate courses (MPC formulation, optimization, stability), the order, the didactic story-line, and the included advanced topics can substantially vary between graduate MPC courses.

\section{Hypotheses and Conclusions} \label{sec:hypo}

\subsection{Hypotheses}
This note has advocated the prospect of teaching practical MPC in undergraduate courses. 
We illustrated that students must meet only limited requirements before they can successfully participate  in an undergraduate MPC course. At the same time, MPC is one of the conceptually
easiest avenues towards MIMO control and towards constrained control. Our main hypotheses are:
\vspace*{0mm}
\begin{itemize}\setlength{\itemsep}{0pt}
\item[(H1)]  Undergraduate MPC teaching is promising and may
increase the attractivity of graduate studies in systems and control. 
\item[(H2)] MPC can be taught immediately after completing basic math courses and a first introductory control course (i.e., in the 2$^{\text{nd}}$ or 3$^{\text{rd}}$ year of a B.Sc. program).\vspace*{0mm}
\end{itemize}

Depending on the curriculum and the focus of the hosting department, (H2) can even be sharpened. For instance, if control is taught to students in computer science, MPC 
could reasonably serve as the first (and possibly only) contact with control. In any case, teaching MPC will also rely on suitable textbooks.  
The brief overview on available MPC and NMPC textbooks---which does not claim to be comprehensive---shows that the community has been very productive. Yet, it stands out that:\vspace*{0mm}
\begin{itemize}
\item[(H3)] There appears to be no dedicated undergraduate text book on MPC, while many recent monographs cover the linear-quadratic material as a gateway to more advanced topics. \vspace*{0mm}
\end{itemize}

In the view of the authors, consensus on which elements to cover in undergraduate MPC course is rather easily obtained. However, teaching of advanced MPC material is a different matter: \vspace*{0mm}
\begin{itemize}\setlength{\itemsep}{0pt}
\item[(H4)] Advanced courses offer synergy potential with numerous topics in systems and control and beyond. If the MPC concept is conveyed in introductory lectures, this provides freedom to capitalize on such synergies. \vspace*{-4mm}
\end{itemize}

\subsection{Final remarks}
We are convinced that MPC as a research field will continue to strive due to its flexibility with respect to outside trends and the inside creativity of the community to explore new directions. Moreover, while the decision how to teach MPC on an advanced graduate level may depend on the specific setting (institution, existing related courses, number and prior knowledge of students, \dots), there appears to be no doubt that transferring elements of MPC teaching to the undergraduate level is a promising, if not necessary, endeavor which may turn into an attractive gateway towards graduate studies in  systems and control.
Likewise, offering MPC education to non-engineering students, e.g., from computer science, is promising as system and control methods are rapidly diffusing into other domains.

\footnotesize
\bibliography{paper}

\end{document}